\documentclass[12pt]{article}
\usepackage{amssymb}
\usepackage{amsthm}
\usepackage{amsmath}
\usepackage{graphicx}

 \newtheorem{thm}{Theorem}[section]
 
 \newtheorem{lem}[thm]{Lemma}
 
 \theoremstyle{definition}
 
 \newtheorem{rem}[thm]{Remark}
 \numberwithin{equation}{section}

\theoremstyle{definition}

\theoremstyle{remark}

\begin{document}
\title{Infinite-time Exponential Growth of the Euler Equation on Two-dimensional Torus}
\author{Zhen
 Lei\footnote{School of Mathematical Sciences; LMNS and Shanghai
 Key Laboratory for Contemporary Applied Mathematics, Fudan University, Shanghai 200433, P. R.China. {\it Email:
 leizhn@gmail.com}}
  \and Jia Shi\footnote{School of Mathematical Sciences, Fudan University, Shanghai 200433, P. R.China. {\it Email: jiashi12@fudan.edu.cn} }
 }
\date{\today}
\maketitle

\begin{abstract}
For any $A > 2$, we construct solutions to the two-dimensional incompressible Euler equations on the torus $\mathbb{T}^2$ whose vorticity gradient $\nabla\omega$ grows exponentially in time: $$\|\nabla\omega(t, \cdot)\|_{L^\infty} \gtrsim e^{At},\quad \forall\ t \geq 0.$$
\end{abstract}

\maketitle





\section{Introduction}

We consider the incompressible Euler equations on the two-dimensional torus
\begin{equation}\label{Euler}
\begin{cases}
u_t + u\cdot\nabla u + \nabla p = 0,\\[-4mm]\\
\nabla\cdot u  = 0.
\end{cases}
\end{equation}
Here $u: [0, \infty)\times \mathbb{T}^2 \longrightarrow \mathbb{R}^2$ is the velocity field and $p: [0, \infty)\times \mathbb{T}^2 \longrightarrow \mathbb{R}$ the scalar pressure. We use $\mathbb{T}^2$ to denote the two-dimensional torus $[-1, 1]^2$. The global regularity of solutions to the two-dimensional Euler equation \eqref{Euler} is known since the work
of Wolibner \cite{Wolibner} and H${ \rm \ddot{o}}$lder \cite{Holder}. More modern and accessible
proofs which hold true for the bounded domain case, the whole space case and also the periodic domain case can be found in \cite{KiselevSverak, MajdaB, Chemin} and the references therein. Recall that classical solutions in Sobolev space $H^{s}$ ($s > 2$) or H${\rm \ddot{o}}$lder space $\mathcal{C}^{1, \gamma}$ ($0 < \gamma < 1$) with initial velocity $u_0 \in L^2$ and $\omega_0 \in L^\infty$ exhibit the following conservation laws
$$\|u(t, \cdot)\|_{L^2} = \|u_0\|_{L^2},\quad \|\omega(t, \cdot)\|_{L^\infty} = \|\omega_0\|_{L^\infty},\quad \forall\ t \geq 0,$$ where $\omega$ denotes the vorticity of the fluid flows. On the other hand, the best known upper bound on the growth of the gradient of vorticity and higher order Sobolev norms is double exponential in time (see for instance, \cite{Yudovich, Wolibner, Holder}).

We are interested in whether there exist solutions of \eqref{Euler} which initially oscillate only on scales comparable to the spatial period and eventually oscillate on arbitrarily short spatial scales for all time $t > 0$. One can quantify such motion in terms of the growth in time of higher Sobolev norms $\|u(t, \cdot)\|_{H^s}$ with $s > 2$, or the $L^\infty$ norm of the vorticity gradient $\|\nabla\omega(t, \cdot)\|_{L^\infty}$. Such kind of problem is motivated by a diverse body of literature on the study of forward energy cascade and weak turbulence theory. See Colliander, Keel, Stanffilani, Takaoka and Tao \cite{Taoetc} for more discussions toward this line of references. In particular, a long standing significant open question is whether such double exponential in time upper bounds for solutions of \eqref{Euler} are sharp (see, for instance, Kiselev and \v{S}ver\'{a}k \cite{KiselevSverak}, Tao \cite{Tao} and also Bourgain \cite{Bourgain}). The main result of this article is the construction of solutions to \eqref{Euler} in the periodic domain case whose vorticity gradient grows exponentially in time for all time $t > 0$. With minor modifications, one can also show similar results for the bounded domain case. See Remark \ref{rem}.

To place our result in context, we review a few highlights from the study on the lower bounds of higher Sobolev norms $\|u(t, \cdot)\|_{H^s}$ with $s > 2$, or the $L^\infty$ norm of the vorticity gradient $\|\nabla\omega(t, \cdot)\|_{L^\infty}$ for the two-dimensional incompressilbe Euler equations. In a series of works \cite{Denisov1, Denisov2, Denisov3}, Denisov began to study this problem. In \cite{Denisov1}, he constructed an example with, in the time average sense, superlinear growth in its vorticity gradient in the periodic case.
In \cite{Denisov2}, he showed that the growth can be double exponential for any given (but finite) period
of time. In \cite{Denisov3}, he constructed a patch solution to 2D Euler equation where, under the action
of a regular prescribed stirring velocity, the distance between the boundaries of two patches
decreases double exponentially in time. We also refer to a discussion at Terry Tao's blog \cite{Tao}
for more information on the problem and related questions.

Very recently, a significant progress on this problem was made by Kiselev and \v{S}ver\'{a}k \cite{KiselevSverak} who constructed smooth solutions on disc with the double exponential growth in the gradient of vorticity in infinite time:
$$\|\nabla\omega(t, \cdot)\|_{L^\infty} \geq \|\omega_0\|_{L^\infty}C_0^{ce^{ct}},\quad {\rm for\ some\ C_0, c\ and}\ \forall\ t \geq 0.$$
We remark that in \cite{KiselevSverak} the boundary is crucial  and the growth happens on the boundary. A formula on velocity field in the small-scale region plays a key role there.  Kiselev and \v{S}ver\'{a}k also pointed out that whether the growth of vorticity gradient for all time is possible in the bulk of the fluid remains a challenging open problem. See also the review of their paper on MathSciNet given by  Secchi (MR3245016). In \cite{AZE}, Zlato\v{s} developed a sharper version of Kiselev-\v{S}ver\'{a}k formula on the velocity field and proved that for any $A \geq 0, \alpha \in (0,1) $, there is $\omega_0\in \mathcal{C}^{1,\alpha}(\mathbb{T}^2)$ with $\|\omega_0\|_{L^\infty} \leq 1$ such that $$\sup_{T_0 \leq t \leq T}\|\omega(t, \cdot)\|_{L^\infty} \geq e^{AT}$$ for all sufficiently large $T_0$.

Using the sharper version of Kiselev-\v{S}ver\'{a}k formula on velocity field given by Zlato${\rm \check{s}}$ \cite{AZE}, we construct solutions to the two-dimensional Euler equation \eqref{Euler} in the periodic domain case with $\mathcal{C}^1$ vorticity whose vorticity gradient has at least exponential growth in time for all $t\geq 0$. Our main result is stated in the following theorem.
\begin{thm}\label{thm}
There exists a constant $K > 0$ such that the following statement is true. For any $A \geq 2$, there is an $\omega_0(x,y)\in \mathcal{C}^1(\mathbb{T}^2$) with $\|\omega_0\|_{L^\infty}\leq 1,\ \|\nabla \omega_0\|_{L^\infty}\leq Ke^{4\sqrt {3}A} $, such that the periodic two-dimensional incompressible Euler equation \eqref{Euler} with initial data $\omega_0$ admits a unique global solution $u(t, \cdot) \in \mathcal{C}^{1, \gamma}$ for any $\gamma \in (0, 1)$ which satisfies $$C_1e^{At/2} \leq \|\nabla\omega(t, \cdot)\|_{L^\infty} \leq \|\nabla\omega_0\|_{L^\infty}\exp\{C_2e^{C_2t}\},\quad \forall\ t \geq 0.$$ Here the positive constant $C_1$ depends only on $A$ and $K$, $C_2$ is given in \eqref{b1} which depends only on the $\mathcal{C}^{1, \gamma}$ norm of the initial velocity.
\end{thm}

\begin{rem}
For smooth bounded domains with two mutually orthogonal axis of symmetry, one can obtain similar exponential growth of the vorticity gradient which happens in the bulk of the fluid. See Remark \ref{rem} for more discussions. At the time of this writting, whether the double exponential growth of vorticity gradient for all time is possible in the bulk of the fluid remains a challenging open problem. In particular, in the whole space case, it seems a challenging problem to construct examples with the exponential growth of vorticity gradient for all time.
\end{rem}

As in \cite{Denisov1, KiselevSverak, AZE}, we will also choose initial data so that the initial vorticity is odd in both $x_1$ and $x_2$. Due to the structure of the Euler equations \eqref{Euler}, this odd property is preserved by vorticity for all later time $t > 0$. Then for any $T > 0$, to show that $\|\nabla\omega(T, \cdot)\|_{L^\infty} \geq C_1e^{\frac{AT}{2}}$, it suffices to show that
$\frac{\omega (T,X(T))}{X_1(T)}$, which is equivalent to $\frac{\omega (0,X(0))}{X_1(T)}$, has the corresponding lower bound $C_1e^{\frac{AT}{2}}$. The trajectory $X(t)$ will be introduced in section 3. For this purpose, we will first show that for the given time $T > 0$, there always exists a material point whose trajectory is sufficiently close to the origin for all time $ t \in [0, T]$. And we will also show that the initial vorticity on this material point is not too small compared with its distance to the $x_2-$axis (at time 0). More precisely, we will show that $\frac{\omega (0,X(0))}{X_1(0)}$  is not too small. On the other hand, we will apply Zlato\v{s}' sharper version of Kiselev-\v{S}ver\'{a}k type formula on velocity field in the periodic domain case to show that $X_1(t)$ decreases in time exponentially for all $0\leq t\leq T$. All of the above  heavily rely on the careful choice of the initial material point $X(0)$ and the construction of the initial data.

\section{Preliminaries}

In this section we first construct our initial data whose velocity field is in $\mathcal{C}^{1, \gamma}$ for any $0 < \gamma < 1$. The initial vorticity gradient has a double logarithmic modulus of continuity. Then we will briefly recall the classical global existence theory of the incompressible Euler equations in $\mathcal{C}^{1, \gamma}$. We will use the $\mathcal{C}^1$ regularity of the vorticity. Since the Euler equations is strongly ill-posed\footnote{The $\mathcal{C}^2$ norm of the velocity field becomes infinite for any $t > 0$ for $\mathcal{C}^2$ initial data.} in $\mathcal{C}^2$ due to the recent important work of Bourgain and Li \cite{B-L}, we will also give an explanation for the $\mathcal{C}^1$ regularity of the vorticity for $t > 0$. At the end of this section we will recall the Kiselev-\v{S}ver\'{a}k type formula of velocity field \cite{KiselevSverak}. Here we use a variant and sharper version given by Zlato\v{s} \cite{AZE}.

\subsection{Set-up of the Initial Data}\label{sub1}

For any given positive constant $A$, we will construct a solution of the two-dimensional incompressible Euler equations so that the vorticity gradient of the solution at least grows exponentially in time: $$\|\nabla\omega(t, \cdot)\|_{L^\infty} \geq C_1e^{At/2},\quad \forall\ t \geq 0.$$ Here $C_1$ is some constant depending only on the initial data and $A$.

Let $\delta \ll 1$ be a positive constant which will be determined later. Let $\omega_0\in \mathcal{C}^1(\mathbb{T}^2)$ be a scalar function which satisfies following constraints:
\begin{itemize}
\item $0\leq \omega_0\leq 1$ in $[0,1]^2$.
\item $\omega_0=1$ in $[\delta,1-\delta]^2$.
\item $\omega_0$ is odd in both $x_1$ and $x_2$.
\item $\omega_0(s,s)=\frac{s}{\log(-\log s)}$ for $s\in[0,\frac{\delta}{2}]$.
\item $\|\nabla \omega_0\|_{L^\infty} \leq \frac{20}{\delta}$.
\end{itemize}
For instance, on $B_\frac{\delta}{2}(0)$ we can choose $\omega_0$ so that $\omega_0(r,\phi)=\frac{r}{\sqrt 2 \log (-\log r)}sin(2\phi)$ in polar coordinates. One may require that $\omega_0 \in \mathcal{C}^\infty(\mathbb{T}^2/\{0\})$.

We remark that the initial vorticity chosen above only possesses the regularity of double logarithmic modulus of continuity. The way of such a choice is partially motivated by the work of Zlato\v{s} \cite{AZE}. It guarantees that the vorticity gradient is still continuous at the origin, which enables us to construct a solution with certain possible growth in the bulk of the fluid.  On the other hand, such a choice of $\omega_0$ also guarantees that the exponential decay of $X_1(t, \alpha)$ (see section 2 for the definition of $X(t, \alpha)$) will eventually lead to the exponential growth of the vorticity gradient. We point out that one may also use $[\ln(-\ln s)]^\alpha$ for any $\alpha > 0$ or $[-\ln(s)]^{\epsilon}$ for some $\epsilon \ll 1$ instead of $\ln(-\ln s)$ in the above choice.

It is easy to check that the initial velocity $u_0 \in \mathcal{C}^{1, \gamma}$ for any $\gamma \in (0, 1)$. By the classical existence theory (see for instance, \cite{MajdaB, KiselevSverak}, etc.), we know that there exists a unique global classical solution $u \in \mathcal{C}([0, \infty), \mathcal{C}^{1, \gamma})$ to the two-dimensional Euler equations with initial data $u_0$.
Moreover, one has
\begin{equation}\label{b1}
\|\nabla u(t, \cdot)\|_{L^\infty} \leq C_2e^{C_2t},
\end{equation}
where $C_2$ is a positive constant depending only on the $\mathcal{C}^{1, \gamma}$ norm of the initial velocity.

\subsection{$\mathcal{C}^1$ Regularity of Vorticity}

We will use the $\mathcal{C}^1$ regularity of the vorticity. Note that our initial velocity may not be a $\mathcal{C}^2$ function. On the other hand, the Euler equations is strongly ill-posed in $\mathcal{C}^2$ due to the recent important work of Bourgain and Li \cite{B-L}. Our observation is that the vorticity is still a $\mathcal{C}^1$ function for all $t > 0$. To see this, we use $X(t, \alpha)$ to denote the flow map determined by the velocity field $u$:
$$\frac{dX(t, \alpha)}{dt} = u(t, X(t, \alpha)),\quad X(0, \alpha) = \alpha.$$
$X^{-1}(t - s, x)$ is the inverse flow map starting from $(t, x)$, which is determined by
$$\frac{dX^{-1}(t - s, x)}{ds} = - u(s, X^{-1}(t - s, x)),\quad X^{-1}(t - s, x)\big|_{s = 0} = x.$$
By \eqref{b1}, it is easy to see that $\nabla X^{-1}$ is a continuous function. Then the $\mathcal{C}^1$ regularity of vorticity is  a straightforward consequence of
\begin{equation}\nonumber
\nabla\omega(t, x) = (\nabla X^{-1}(t, x))^T\nabla\omega_0(X^{-1}(t, x))
\end{equation}
which can be derived from the following formula for vorticity gradient:
\begin{equation}\label{equ}
\partial_t(\nabla\omega) + u\cdot\nabla(\nabla\omega) = - (\nabla\omega)\cdot\nabla u
\end{equation}

\subsection{The Kiselev-\v{S}ver\'{a}k Type Formula of Velocity}

To get the the exponential decay in time of $X_1(t, \alpha)$, we need the Kiselev-\v{S}ver\'{a}k type formula of velocity field. Here we choose to use a sharper version on the torus given by Zlato\v{s} \cite{AZE}. In the lemma below, $Q(x)\ :=[x_1,1]\times[x_2,1]$.
\begin{lem}[Zlata\v{s}]\label{lem1}
Let $\omega_0 \in L^\infty(\mathbb{T}^2)$ be odd in both $x_1$ and $x_2$. If $x_1$, $x_2$ $\in [0,\frac{1}{2})$, then
\begin{equation}\label{b2}
u_j(t,x)=(-1)^j\Big(\frac{4}{\pi}\int_{Q(2x)}\frac{y_1 y_2}{|y|^4}\omega(t,y)dy+B_j(t,x)\Big)x_j,\quad {\rm for}\ j=1, 2,
\end{equation}
where, with some universal $C_3$, $B_1(t, x)$ and $B_2(t, x)$ satisfy
\begin{eqnarray}\nonumber
|B_1| &\leq& C_3\|\omega_0\|_{L^\infty}\Big(1 + \min{\Big\{\log(1+\frac{x_2}{x_1}),
 x_2\frac{\|\nabla \omega(t,\cdot)\|_{L^\infty}}{\|\omega_0\|_{L^\infty}}\Big\}}\Big),\\\nonumber
|B_2| &\leq& C_3\|\omega_0\|_{L^\infty}\Big(1+\min{\Big\{\log(1+\frac{x_1}{x_2}),x_1\frac{\|\nabla \omega(t,\cdot)\|_{L^\infty}}{\|\omega_0\|_{L^\infty}}\Big\}}\Big).
\end{eqnarray}
\end{lem}

\section{Proof of Theorem \ref{thm}}

The lower bound of the first term in the large bracket of \eqref{b2} is a crucial observation in the work of Kiselev and \v{S}ver\'{a}k in \cite{KiselevSverak}. The following estimate is motivated by \cite{KiselevSverak} and is also very crucial in this article. Now we take $\omega_0$ as the one given in subsection \ref{sub1} Note that $\omega_0\geq 0$ on $[0,1]^2$, and $\omega_0\leq 1$ on a subset of $[0,1]^2$ with measure less than $4\delta$. By the incompressibility, one sees that $\omega(t,x)$ has such same properties as $\omega_0$. Hence, for $x\in[0,\delta_1/2]^2$ ($\delta_1 \leq \frac{1}{4}$ to be chosen later), we have
\begin{eqnarray}\nonumber
\frac{4}{\pi}\int_{Q(2x)}\frac{y_1 y_2}{|y|^4}\omega(t,y)dy &\geq&
\frac{4}{\pi}\int_{2\delta_1}^1\int_\frac{\pi}{6}^{\frac{\pi}{3}}\frac{1}{r}\sin\theta\cos\theta\omega(r,\theta)dr d\theta\\\nonumber
&\geq&\frac{\pi}{6}\frac{\sqrt 3}{\pi}\int_{\sqrt{4\delta_1^2+\frac{48}{\pi}\delta}}^1\frac{1}{r}dr\\\nonumber
&=&\frac{\sqrt 3}{12}[-\ln(4{\delta_1}^2+\frac{48}{\pi}\delta)]
\end{eqnarray}
Now we take $\delta = \frac{\pi}{96}e^{-4\sqrt 3(A+2C_3)}$ and $\delta_1 = \min{\{\frac{\sqrt 2}{4}e^{-2\sqrt3 (A+2C_3)},\frac{\delta}{2}\}} $ to derive that
\begin{equation}\label{b3}
\frac{4}{\pi}\int_{Q(2x)}\frac{y_1 y_2}{|y|^4}\omega(t,y)dy\geq (A+2C_3).
\end{equation}
The constant $K$ in the Theorem \ref{thm} can just be taken as $K = \frac{1920}{\pi}e^{8\sqrt {3}C_3}$. We will show that $|B_1(t, X(t))|$ is controlled by $2C_3$ for a careful chosen trajectory $X(t)$ below.

Now for any large $T$, we need to show that $$\|\nabla\omega(T, \cdot)\|_{L^\infty} \geq C_1e^{AT/2}.$$ Our strategy is to find an initial material point and track its motion, for this given time $T$. Let $s > 0$ be a small constant depending on $T$ and being determined below. We study the trajectory $X(t) = X(t, (s, s))$ which starts from the material point $X(0, (s, s))=(s,s)$. For all $0\leq t \leq T$, we claim that $$X_2(t)\cdot\|\nabla \omega(t, \cdot)\|_{L^\infty}\leq 1,$$ and $$X(t)\in [0,\delta_1]^2,$$
if the initial material point $(s, s)$ is appropriately chosen.

Indeed, we first note that
\begin{equation}\nonumber
\frac{dX_1(t)}{dt} = u_1(t, X(t)),\quad
\frac{dX_2(t)}{dt} = u_2(t, X(t)).
\end{equation}
Using the facts that $u_1(t, 0, X_2(t)) = 0$ and $u_2(t, X_1(t), 0) = 0$, and the estimate for $\nabla u$ in \eqref{b1}, one has
\begin{eqnarray}\nonumber
|\frac{dX_1}{dt}| &\leq& X_1 \cdot \|\nabla u\|_{L^\infty}\leq C_2 X_1e^{C_2t},\\\nonumber
|\frac{dX_2}{dt}| &\leq& X_2 \cdot \|\nabla u\|_{L^\infty}\leq C_2X_1e^{C_2t}.
\end{eqnarray}
Then Gronwall's inequality gives that
\begin{eqnarray}\label{X}
 X_1(t) &\leq& s\cdot \exp{\{e^{C_2 t}-1\}}\leq s\cdot\exp{\{e^{C_2 T}-1\}},\\\nonumber
 X_2(t) &\leq& s\cdot \exp{\{e^{C_2 t}-1\}}\leq s\cdot \exp{\{e^{C_2 T}-1\}}.
\end{eqnarray}
Moreover, it is easy to derive from \eqref{equ} and Gronwall's inequality that $$\|\nabla\omega(t, \cdot)\|_{L^\infty} \leq \|\nabla\omega_0\|_{L^\infty}e^{\int_0^t\|\nabla u(s, \cdot)\|_{L^\infty}ds}.$$ Hence, we have
\begin{eqnarray}\nonumber
X_2(t)\|\nabla \omega(t,\cdot)\|_{L^\infty} &\leq& (s\cdot\exp{\{e^{C_2 t}-1\}})\cdot(
\|\nabla\omega_0\|_{L^\infty}\exp{\{e^{C_2 t}-1\}})\\\nonumber
&\leq& s\cdot \|\nabla\omega_0\|_{L^\infty}\cdot \exp{\{2e^{C_2 T}-2\}}\\\nonumber
&\leq& s\cdot\frac{20}{\delta}\exp{\{2e^{C_2 T}-2\}}.
\end{eqnarray}
Now we choose $s$ as follows
\begin{equation}\label{s}
s=\min{\big\{\frac{\delta}{20}\exp(2-2e^{C_2T}),\quad \frac{\delta_1}{2} \exp(1-e^{C_2T})\big\}}.
\end{equation}
It is clear that such a choice for $s$ guarantees that $X(t)\in[0,\delta_1/2]^2$ in view of \eqref{X}. Moreover, one also has
\begin{eqnarray}\nonumber
X_2(t)\|\nabla \omega(t,\cdot)\|_{L^\infty} \leq 1,\quad \forall 0 \leq t \leq T.
\end{eqnarray}
We have proved our claim.

Now using the above claim, and noting $\|\omega_0\|_{L^\infty} = 1$, one has
$$X_2(t)\frac{\|\nabla \omega(t,\cdot)\|_{L^\infty}}{\|\omega_0\|_{L^\infty}} \leq 1,\quad \forall\ 0 \leq t \leq T.$$
Taking $(t, x) = (t, X(t))$ in Lemma \ref{lem1}, one has $$|B_1(t,X(t))|\leq 2C_3$$ for all $0\leq t\leq T$.
Consequently, by \eqref{b2}, we have
\begin{equation}\nonumber
- u_1(t, X(t)) \geq \Big(\frac{4}{\pi}\int_{Q(2x)}\frac{y_1 y_2}{|y|^4}\omega(t,y)dy - 2C_3\Big)X_1(t),\quad \forall\ 0 \leq t \leq T.
\end{equation}
Recall the estimate in \eqref{b3}, we finally arrive at
\begin{eqnarray}\nonumber
\frac{d}{dt}X_1(t) = u_1(t,X(t)) \leq   -AX_1(t),\quad \forall\ 0\leq t\leq T.
\end{eqnarray}
which gives that
\begin{eqnarray}\label{b6}
X_1(t) \leq se^{-At},\quad \forall\ 0\leq t\leq T.
\end{eqnarray}
This decay in time estimate of $X_1(t)$ and the double logarithmic modulus of continuity of $\nabla \omega_0$ are crucial for the proof of theorem \ref{thm} and eventually lead to the exponentially growth of $\|\nabla \omega(t,\cdot)\|_{L^\infty}$.

Indeed, since $\omega(t,x)$ is odd in $x_1$, one has $\omega(0,x_2,t)=0$. Hence, using \eqref{b6}, we can get
\begin{eqnarray}\nonumber
\|\nabla \omega(T, \cdot)\|_{L^\infty} &\geq& \frac{\omega(T, X(T)) - \omega(T, 0, X_2(t))}{X_1(T)}\\\nonumber
&=& \frac{\omega_0(s, s)}{X_1(T)} \geq \frac{\omega_0(s, s)}{se^{-AT}}\\\nonumber
&=& \frac{e^{AT}}{\log(-\log(s))}.
\end{eqnarray}
According to \eqref{s}, one has
\begin{eqnarray}\nonumber
\log(-\log(s)) &=& \max\big\{\ln\big(\ln\frac{20}{\delta} + 2e^{C_2T} - 2\big),\quad \ln\big(\ln\frac{2}{\delta_1} + e^{C_2T} - 1\big\}\\\nonumber
&\leq& C_4(1 + T)
\end{eqnarray}
for some constant $C_4$ depending only on initial data. Then we have
\begin{eqnarray}\nonumber
\|\nabla \omega(T, \cdot)\|_{L^\infty} \geq \frac{e^{AT}}{C_4(1 + T)}.
\end{eqnarray}
So the proof of Theorem \ref{thm} is finished by choosing $C_1 = \frac{1}{2C_4}$.

\begin{rem}\label{rem}
For the regular bounded domain $\Omega$ with two mutually orthogonal axis of symmetry, we can use the decomposition of the Green function $G_{D}(x,y) = - \frac{1}{2\pi}\log|x-y| + g(x,y)$ to get the estimate of the velocity field $u$ as in Lemma \ref{lem1}. Here $g$ is a harmonic function. Then with minor suitable modifications, the method of this paper can be used to get the similar results for the two-dimensional incompressible Euler equation on certain bounded domains. Again, the growth of the vorticity gradient happens in the bulk of the fluid.
\end{rem}

\section*{Acknowledgement.}

The authors were in part supported by NSFC (grant No. 11421061 and 11222107), National Support Program for Young Top-Notch Talents, Shanghai Shu Guang project, and SGST 09DZ2272900.

\end{document}